\documentclass[12pt]{amsart}

\usepackage{graphicx}
\usepackage{amsthm}
\swapnumbers

\newtheorem{Thm}{Theorem}
\newtheorem{Def}[Thm]{Definition}
\newtheorem{Coro}[Thm]{Corollary}
\newtheorem{Lem}[Thm]{Lemma}
\newtheorem{Prop}[Thm]{Proposition}
\newtheorem{Que}[Thm]{Question}

\begin{document}

\title{Bridge Number and the Curve Complex}
\author{Jesse Johnson}
\address{\hskip-\parindent
        Mathematics Department\\
        University of California\\
        Davis, CA 95616\\
        USA}
\email{jjohnson@math.ucdavis.edu}

\thanks{Research supported by NSF VIGRE grant 0135345}

\begin{abstract}
We show that there are hyperbolic tunnel-number one knots with arbitrarily high bridge number and that ``most" tunnel-number one knots are not one-bridge with respect to an unknotted torus.  The proof relies on a connection between bridge number and a certain distance in the curve complex of a genus-two surface.
\end{abstract}

\maketitle

\section{Introduction}
\label{introsect}

Let $S$ be a sphere smoothly embedded in $S^3$ and let $K$ be a knot transverse to $S$.  The complement in $S^3$ of $S$ consists of two open balls, $B_1$ and $B_2$.  The knot $K$ is in \textit{bridge position} (with respect to $S$) if $K \cap B_1$ is properly isotopic into $S$ and $K \cap B_2$ is properly isotopic into $S$.  The \textit{bridge number} $b_0(K)$ of $K$ is the smallest possible number of components of $K' \cap B_1$ where $K'$ is in bridge position and $K'$ is isotopic to $K$.

The bridge number of the unknot is 0 because the unknot can be isotoped into $B_2$ so that it is parallel to $S$.  We will say that a knot is non-trivial if it is not the unknot.

It is known that there are knots with arbitrarily large bridge number.  In particular, Schubert showed that the bridge number of a satellite knot is bounded below in terms of the bridge numbers of the knots used to construct it.  This fact can be used to construct examples of knots with high bridge number, all containing essential tori or annuli in their complements.

One can also consider the bridge number $b_1(K)$ of a knot $K$ with respect to an
unknotted torus in $S^3$.  Moriah and Rubinstein~\cite{mr:neg}, Morimoto, Sakuma and Yokota~\cite{MSY} and Eudave-Munoz~\cite{eud:tor}  have shown that there are tunnel-number one knots which are not one-bridge with respect to an unknotted torus.  However, there are no known examples of tunnel-number one knots with toroidal bridge number greater than two.  

We will show that tunnel-number one knots with large bridge number and with $b_1(K) > 1$ are fairly common.  In particular, we have the following Theorem:

\begin{Thm}
\label{mainthm}
For every integer $N$, there is a hyperbolic, tunnel-number one knot such that $b_0(K) > N$ and $b_1(K) > 1$.
\end{Thm}

To prove Theorem~\ref{mainthm}, we define a measure of complexity $d(\tau)$ for an unknotting tunnel $\tau$ of a knot $K$.  The main theorem then follows from three lemmas:

\begin{Lem}
\label{bridgelem}
For a non-trivial knot $K$ with an unknotting tunnel $\tau$, $b_0(K) \geq d(\tau)$.
\end{Lem}

\begin{Lem}
\label{torbridgelem}
If $\tau$ is an unknotting tunnel for a knot $K$ and $d(\tau) > 5$ then $b_1(K) > 1$.
\end{Lem}

\begin{Lem}
\label{highNlem}
For every $N$ there is a tunnel-number one knot with unknotting tunnel $\tau$ such that $d(\tau) > N$.
\end{Lem}

The integer $d(\tau)$ is the distance from a certain loop in a handlebody set in the curve complex to a second handlebody set.  Because a handlebody set is fairly homogeneous, Lemma~\ref{highNlem} implies that ``most'' unknotting tunnels have high distance, so ``most'' knots are not one-bridge with respect to an unknotted torus.

By the following theorem of Morimoto, Lemma~\ref{torbridgelem} also tells us about the behavior of these knots under connect sum:

\begin{Thm}[Morimoto~\cite{Mor}]
Let $K_1$ and $K_2$ be knots in $S^3$ with tunnel-number one.  Then the tunnel number of their connect sum is 3 if and only if $b_1(K_1) > 1$ and $b_1(K_2) > 1$.
\end{Thm}

\begin{Coro}
If knots $K_1$ and $K_2$ have unknotting tunnels $\tau_1$ and $\tau_2$, respectively, such that $d(\tau_1) > 5$ and $d(\tau_2) > 5$ then the connect sum of $K_1$ and $K_2$ has tunnel-number 3.
\end{Coro}

We define the complexity $d(\tau)$ in Section~\ref{defsect} and describe some of its properties in Section~\ref{stdistsect}.  We prove Lemma~\ref{bridgelem} in Section~\ref{bridgesect} and Lemma~\ref{torbridgelem} is proved in Section~\ref{torbridgesect}.  In Section~\ref{ccompsect}, we discuss the geometry of the curve complex, and in Section~\ref{mainthmsect} we prove Lemma~\ref{highNlem} and Theorem~\ref{mainthm}.

Thanks to David Futter for a number of useful suggestions.

\section{Unknotting Tunnels and the Curve Complex}
\label{defsect}

A \textit{handlebody} is a manifold with boundary that is homeomorphic to a regular neighborhood of a connected graph embedded in a closed, orientable 3-manifold.  The genus of a handlebody is the genus of its boundary.  Let $H$ be a genus-$g$ handlebody and $\Sigma$ a closed, orientable, genus-$g$ surface.  Let $\phi : \Sigma \rightarrow \partial H$ be a homeomorphism.

\begin{Def}\textup{
The \textit{curve complex} $C(\Sigma)$ is the graph whose vertices are isotopy
classes of simple closed curves in $\Sigma$ and edges connect vertices
corresponding to disjoint curves.
}\end{Def}

For more detailed descriptions of the curve complex, see~\cite{hemp:cc} and~\cite{mm:hyp}.

\begin{Def}\textup{
The \textit{handlebody set} $\mathbf{H} \subset C(\Sigma)$ corresponding to $(H,\phi)$ is the set of vertices $\{l \in C(\Sigma) : \phi(l)$ bounds a disk in $H\}$.
}\end{Def}

We will not mention the map $\phi$ when it is obvious from the context. Given vertices $l_1, l_2$ in $C(\Sigma)$, the distance $d(l_1, l_2)$ is the geodesic distance: the number of edges in the shortest path from $l_1$ to $l_2$.  This definition  extends to a definition of distances between subsets $A,B$ of $C(\Sigma)$ by defining $d(A,B) = \min \{d(a,b) : a \in A, b \in B\}$ and for distances between a point and a set similarly.

Given a knot $K$, an \textit{unknotting tunnel} is an arc $\tau$ with its endpoints in $K$ such that for a regular neighborhood $N$ of $K \cup \tau$, $S^3 \setminus N$ is a genus-two handlebody $H_2$.  We call $K$ a \textit{tunnel-number one} knot if there is an unknotting tunnel for $K$.  The graph $K \cup \tau$ defines a Heegaard splitting $(\Sigma, H_1, H_2)$ of $S^3$ where $\Sigma = \partial H_2$ and $H_1 = N \cup \Sigma$.  For $i = 1,2$, let $\mathbf{H_i} \subset C(\Sigma)$ be the handlebody set induced by $H_i$ with the inclusion map $\Sigma \rightarrow \partial H_i$.  

Let $D_1, D_2, D_3 \subset H_1$ be the meridian disks dual to the three edges of $K \cup \tau$ and assume $D_1$ is the disk dual to $\tau$.  In other words, the disks should be properly embedded, pairwise disjoint and each disk should intersect the corresponding edge in a single point and the other edges in no points.  

\begin{figure}[htb]
  \begin{center}
  \includegraphics[width=3.5in]{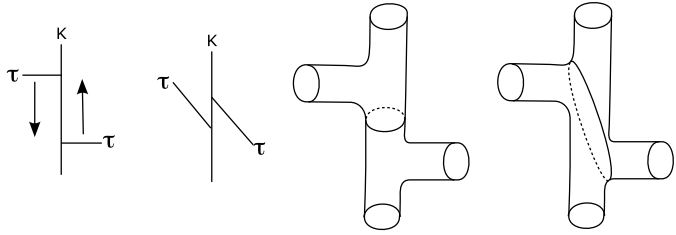}
  \caption{Sliding the ends of $\tau$ past each other corresponds to replacing the disk $D_2$ or $D_3$.}
  \label{slidefig}
  \end{center}
\end{figure}

\begin{Def}\textup{
A \textit{tunnel isotopy} of $\tau$ is a sequence of isotopies of $K \cup \tau$ that fix $K$ and moves in which the ends of $\tau$ are slid past each other along $e_1$ or $e_2$, as in Figure~\ref{slidefig}.  
}\end{Def}

After the ends of $\tau$ are slid past each other, a set of meridian disks for the new graph can be found by replacing $D_2$ or $D_3$ with a new disk.

\begin{Lem}
\label{tunslidelem}
Given two unknotting tunnels, the induced Heegaard splittings of the knot complement are isotopic if and only if the tunnels are slide isotopic.
\end{Lem}

The proof of the Lemma will be left to the reader.  Define $d(\tau)$ to be the distance $d(\partial D_1, \mathbf{H_2})$ in $C(\Sigma)$.  The slide moves defined above do not change the isotopy class of the disk $D_1$ dual to $\tau$ so $d(\tau)$ is invariant under tunnel isotopies.

\section{Standard distance}
\label{stdistsect}

For a Heegaard splitting $(\Sigma, H_1, H_2)$ of an arbitrary 3-manifold $M$, the \textit{standard distance} is $d(\Sigma) = d(\mathbf{H_1}, \mathbf{H_2})$, where $\mathbf{H}_1$ and $\mathbf{H}_2$ are the handlebody sets corresponding to $H_1$ and $H_2$, respectively.  

For the genus-two Heegaard splitting $(\Sigma, H_1, H_2)$ of $S^3$ defined by $K \cup \tau$, the complement in $H_1$ of a regular neighborhood $N$ of $K$ is a compression body.  Thus the inclusion into $M \setminus N$ of $(\Sigma, (H_1 \setminus N), H_2)$ is a genus-two Heegaard splitting of the knot complement and we can consider the standard distance $d(\Sigma)$ of this splitting.

The meridian disk $D_1$ dual to $\tau$ is a properly embedded, essential disk in $H_1 \setminus N$ so $\partial D_1$ defines an element of the handlebody set $\mathbf{H}_1$.  Thus we immediately have the inequality $d(\Sigma) \leq d(\tau)$.  The reverse inequality does not necessarily hold, but it isn't far off.

\begin{Lem}
\label{splitlem}
If $\Sigma$ is the Heegaard splitting of the complement of $K$ induced by an unknotting tunnel $\tau$ then $d(\Sigma) \geq d(\tau) - 1$.
\end{Lem}

\begin{proof}
Let $D_1$ be the meridian disk for $\tau$ in $H_1$ and let $D$ be any properly embedded, essential disk in $H_1$. We will show that $D$ can be isotoped disjoint from $D_1$, implying $d(\partial D, \partial D_1) \leq 1$ and $d(\partial D, \mathbf{H_2}) \geq d(\partial D_1, \mathbf{H_2}) - 1$.

Isotope $D$ so as to minimize the number of components of $D \cap D_1$. Since $H_1$ is irreducible, we can eliminate any loop components, implying that $D \cap D_1$ consists of a number of arcs.  If the intersection is empty then we're done.  Otherwise, consider an outermost arc, cutting off a disk $D' \subset D$.  This disk lies in $H_1 \setminus D_1$, which is homeomorphic to $T^2 \times [0,1]$, where $T^2$ is a torus.  Every disk in this manifold is boundary parallel. 

Since $D'$ is boundary parallel, the loop $\partial D'$ bounds a disk in some component of $T^2 \times [0,1]$.  Let $D'' \subset T^2 \times \{0\}$ be this disk.  The image in $T^2 \times \{0\}$ of $D_1$ is a pair of disks.  One of these disks intersects $\partial D'$ in an arc and the other is disjoint from $\partial D'$.  Let $D'_1$ be the disk disjoint from $\partial D'$ and let $D''_1$ be the second disk.  Then $D'_1$ is either contained in $D''$ or disjoint from $D''$.

If $D_1'$ is not contained in $D''$ then the arc of $\partial D'$ that is disjoint from $D_1$ can be slid across $D''$.  The induced isotopy of $D$ reduces the number of intersections with $D_1$.

Assume for contradiction $D'_1$ is contained in $D''$.  Let $A$ be the annulus $D'' \setminus D'_1$.  Then $\partial D \cap A$ consists of a number of arcs.  Any trivial arc in $A$ can be isotoped into $D'_1$ or $D''_1$.  Since $D \cap D_1$ was minimized, every arc of $D \cap A$ must be essential.  In particular, every arc must go from $D'_1$ to $D''_1$.  

There is also an arc of  $\partial D'$ which goes from $D''_1$ back to itself, so $\partial D$ intersects $D''_1$ in more points than $D'_1$.  This is impossible because both disks are images of the same disk $D_1$.  This contradiction implies that $D$ and $D_1$ must be disjoint.
\end{proof}

This Lemma allows us to immediately employ a number of useful results about standard distance to study unknotting tunnels.

\begin{Thm}[Scharlemann and Tomova~\cite{tom:dist}]
\label{sctmthm}
If $(\Sigma, H_1, H_2)$ and $(\Sigma', H'_1, H'_2)$ are genus $g$ Heegaard splittings of $M$ and $d(\Sigma) > 2g$ then $\Sigma'$ is isotopic to $\Sigma$.
\end{Thm}

Combining this with Lemma~\ref{splitlem} and Lemma~\ref{tunslidelem}, we conclude the following:

\begin{Coro}
\label{sctmcoro}
If $d(\tau) > 5$ then any unknotting tunnel $\tau'$ for $K$ is slide isotopic to $\tau$.
\end{Coro}


We can also find a nice connection between $d(\tau)$ and the Seifert genus of $K$.

\begin{Prop}[Scharlemann and Thompson~\cite{thm:rho}]
\label{sgen1lem}
Suppose $K$ is a knot and $\tau$ an unknotting tunnel.   Then $\tau$ can be made disjoint from a minimum-genus Seifert surface by a slide isotopy.
\end{Prop}

The following Lemma is a special case of the main Theorem of \cite{sch:first}.

\begin{Lem}[Scharlemann~\cite{sch:first}]
\label{sgen2lem}
Let $l \subset \Sigma$ be a simple closed curve which bounds a properly embedded surface $S$ in $H_2$ then $d(l,\mathbf{H_2}) \leq 1-\chi(S)$.
\end{Lem}

\begin{Coro}
If $h$ is the Seifert genus of $K$ then $d(\tau) \leq 2h + 1$.
\end{Coro}

\begin{proof}
By Proposition~\ref{sgen1lem}, take $\tau$ to be disjoint from some Seifert surface $S$.  Then $S \cap H_1$ is an annular neighborhood of $\partial S$ and $S' = S \cap H_2$ is a properly embedded, once-punctured, genus $h$ surface in $H_2$.  Moreover, $\partial D_1$ is disjoint from $l = \partial S'$ so $d(l, \partial S') \leq 1$.  By Lemma~\ref{sgen2lem}, $d(l,\mathbf{H}_2) \leq 1-\chi(S') = 2h$ so $d(\partial D_1, \mathbf{H_2}) \leq 2h+1$.
\end{proof}

\section{Bridge number}
\label{bridgesect}

If $K$ is the unknot then its complement is a solid torus.  Let $D$ be an embedded disk in $S^3$ whose boundary is $K$ and let $\alpha$ be an arc in $D$ whose endpoints are in $K$.  Let $\tau$ be the result of isotoping the interior of $\alpha$ away from $D$.   This $\tau$ is an unknotting tunnel for the unknot because the complement of an open neighborhood $K \cup \alpha$ is a genus-two handlebody so $\alpha$ is an unknotting tunnel and $\tau$ is isotopic to $\alpha$.  The disk $D$ is properly embedded in the complement and is disjoint from a meridian of $\tau$ so $d(\tau) = 1$.

By Lemma 2.7 in \cite{sct:hbody}, any two Heegaard splittings of a handlebody are isotopic. Thus Lemma~\ref{tunslidelem} implies that any unknotting tunnel $\tau'$ for the unknot is slide isotopic to $\tau$ and $d(\tau') = 1$.  (One could also Theorem 1' in~\cite{grdn}, which states that any arc in a handlebody whose complement is a handlebody is a boundary parallel arc.)

The converse is also true:  If $d(\tau) = 1$ then there is a properly embedded, essential disk $D$ in the complement of a neighborhood $N$ of $K \cup \tau$ whose boundary is disjoint from the meridian $D_1$ of $\tau$.  Then $D$ is properly embedded in the complement of $N \setminus D_1$.  The set $N \setminus D_1$ is a regular neighborhood of $K$ so the existence of the disk $D$ implies that $K$ is the unknot.  So, $K$ is the unknot if and only if $d(\tau) = 1$.

Before continuing with the proof of Lemma~\ref{bridgelem}, we need to state one more lemma.  Let $K$ be a knot in minimal bridge position with respect to a sphere $S$ and let $\tau$ be an unknotting tunnel for $K$.

\begin{Lem}[Goda, Scharlemann and Thompson~\cite{gst:level}]
\label{gstlem}
There is a tunnel isotopy taking $\tau$ into $S$, if we allow the last move to slide the endpoints of $\tau$ together in $K$.
\end{Lem}

We will use this in the following proof.

\begin{proof}[Proof of Lemma~\ref{bridgelem}]
Induct on the bridge number $n = b_0(K)$.  If $K$ has bridge number less than $2$ then $K$ is the unknot and $d(\tau) = 1$.
 
For $n \geq 2$, assume the result is true for every knot with bridge number strictly less than $n$.  Let $K$ be a knot with bridge number $b_0(K) = n$, with unknotting tunnel $\tau$.  By Lemma~\ref{gstlem}, $\tau$ can be made level by a tunnel isotopy if we allow that the last move takes both endpoints of the tunnel to the same point in $K$.  Assume that $\tau$ has been leveled in this fashion.

If the endpoints of $\tau$ are the same point in $K$ then $\tau$ forms a loop in the sphere $S$ so $\tau$ is the unknot.   The loop $K$ can be thought of as an edge with both ends attached to the same point on the loop $\tau$. Sliding the ends of this edge away from each other along $\tau$ turns $K$ into an arc $\tau'$ which is an unknotting tunnel for $\tau$.  Since $\tau$ is the unknot, $d(\tau') = 1$.  The meridian disk for $\tau'$, is disjoint from the meridian of $\tau$ (both in $H_1$) so $2 \geq d(\tau)$.  Since we assumed $b_0(K) \geq 2$, we have $b_0(K) \geq d(\tau)$.

Assume $\tau$ sits in a level sphere and has distinct endpoints.  Then $\tau \cap K$ consists of two points which separate $K$.  Let $e_2$ and $e_3$ be the (arc) components of $K \setminus \tau$. (These correspond to the meridian disks $D_2$ and $D_3$.)  Then $K' = \tau \cup e_2$ and $K'' = \tau \cup e_3$ are loops.

By isotoping $\tau$ slightly we can put the loop $K'$ in bridge position.  A different isotopy of $\tau$ will put $K''$ in bridge position.  The number of bridges in $e_2$ and $e_3$ add up to the number of bridges in $K$, so one of the loops $K'$ and $K''$ will have bridge number strictly less than $K$.  Without loss of generality, assume $b_0(K') < b_0(K)$.

The edge $e_3$ forms an unknotting tunnel for $K'$.  The meridian disks $D_1, D_3$ are disjoint so the vertices in $C(\Sigma)$ corresponding to their boundaries are connected by an edge and $d(\tau) \leq d(e_3) +1$.  Because $b_0(K') < n$, we conclude that $b_0(K') \geq d(e_3)$ or $d(e_3) = 1$.   In either case, $b_0(K) \geq b_0(K') + 1 \geq d(e_3) + 1 \geq d(\tau)$.
\end{proof}

\section{Toroidal Bridge Number}
\label{torbridgesect}

The method used to prove Lemma~\ref{bridgelem} cannot be generalized to toroidal bridge number because it is not known if unknotting tunnels can always be leveled with respect to an unknotted torus.  For knots with $b_1(K) = 1$, however, we can get around this problem using the following Lemma:

\begin{Lem}
\label{leveltunlem}
If $K$ is a knot such that $b_1(K) = 1$ then $K$ has an unknotting tunnel, $\tau'$, which sits in a level torus and such that $d(\tau') \leq 2$.
\end{Lem}

\begin{proof}
Let $K$ be a knot in one-bridge position with respect to an unknotted torus $T$.  Let $B_0, B_1 \subset S^3$ be the solid tori defined by $T$.  Let $N'$ be a regular neighborhood of $K$.  There is a disk $D \subset B_1$ whose boundary consists of the arc $K \cap B_1$ and an arc $\beta$ in $T$.  (This is because $K$ is in bridge position with respect to $T$.)  A regular neighborhood $N$ of $N' \cup B_0 \cup D_1$ is isotopic to $B_0$ so its complement is a solid torus.  By removing the disk $D_1$ from the union and attaching it to the complement, we find that the complement of $N' \cup B_0$ is a genus-two handlebody.

Let $D_0$ be a properly embedded essential disk in $B_0$ that is disjoint from $N'$ and let $D_1$ be a similar disk in $B_1$. Such disks exist because $K$ is in bridge position.  Let $D_2$ be a meridian of $N'$ which is disjoint from $B_0$.  Because $D_0$ and $D_1$ are disjoint from $N'$, they are properly embedded in $N$ and in the complement of $N$, respectively.  The disk $D_2$ is also properly embedded and essential in $N$.  The boundary of each disk is disjoint from $D_2$ in $\partial N$ so $d(\partial D_0, \partial D_1) \leq 2$.

The complement in $N$ of an open neighborhood of $D_0$ is a solid torus which is isotopic to a regular neighborhood of $K$.  Thus there is an arc $\tau$ in $N$ which is dual to $D_0$ and such that $\tau$ is an unknotting tunnel for $K$. We then have $d(\tau) = d(\partial D_0, \mathbf{H}_2) \leq d(\partial D_0, \partial D_1) \leq 2$.
\end{proof}

Note that the complement in $N$ of the disk $D_2$ is isotopic to $B_0$.  By sliding the endpoints of $\tau$ together along the arc $K \cap B_0$, we can turn $\tau$ into a loop which is a core of $B_0$ and therefore the unknot.  In other words, if $K$ is one-bridge with respect to the torus then it has an unknotting tunnel which is an unknotted loop.  The converse of this statement is true as well, but is not necessary for this paper, so the proof of the following Theorem is left to the reader:

\begin{Thm}
A knot $K$ is one-bridge with respect to an unknotted torus if and only if $K$ has an unknotting tunnel which is an unknotted loop.
\end{Thm}

The existence of a level tunnel, alone, is not good enough to prove a general bound on the bridge number.  We need to know that an arbitrary tunnel for the knot has bounded distance.  Luckily, if $d(\tau)$ is high enough, Corollary~\ref{sctmcoro} implies that any unknotting tunnel is isotopic to $\tau$.

\begin{proof}[Proof of Lemma~\ref{torbridgelem}]
Let $K$ be a knot with unknotting tunnel $\tau$ such that $d(\tau) > 5$.  Assume for contradiction $b_1(K) = 1$.  By Lemma~\ref{leveltunlem}, there is an unknotting tunnel $\tau'$ such that $d(\tau') \leq 2$.  However, because $d(\tau) > 5$, Corollary~\ref{sctmcoro} implies $\tau'$ is slide isotopic $\tau$ and this contradicts the assumption that $d(\tau) > 5$ while $d(\tau') \leq 2$.  The contradiction completes the proof.
\end{proof}

Considering Lemma~\ref{torbridgelem} along with the proof of Lemma~\ref{bridgelem} raises the following question:

\begin{Que}
\label{bigquestion}
Can unknotting tunnels be leveled with respect to an unknotted torus? 
\end{Que}

A positive answer would allow us to replace $b_0(\tau)$ with $b_1(\tau)$ in the upcoming proof of Theorem~\ref{mainthm}, implying the following:

\begin{Prop}
If the answer to Question~\ref{bigquestion} is ``yes'' then there are tunnel-number-one knots with arbitrarily high bridge number with respect to an unknotted torus.
\end{Prop}

\section{Geometry of the Curve Complex}
\label{ccompsect}

To find an appropriate knot and unknotting tunnel, we will first find a disk $D$ in a standardly embedded handlebody $H_1 \subset S^3$ such that $d(\partial D, \mathbf{H_2}) > N$.  From $D$ we will construct a knot $K$ with unknotting tunnel $\tau$ whose meridian disk is $D$.  In this section, we will prove the following Lemma:

\begin{Lem}
\label{unbndlem}
For every $N \in \mathbf{N}$, there is a vertex $v \in \mathbf{H_1}$ such
that $d(v, \mathbf{H_2}) > N$.
\end{Lem}

To prove this we will employ a number of powerful results about the geometry of the curve complex.

\begin{Def}\textup{
Given sets $X$ and $Y$ in a metric space, we will say that \textit{$X$ is an unbounded distance from $Y$} if for every $N$, there is an $x \in X$ such that $d(x, Y) > N$.  Otherwise, \textit{$X$ is a bounded distance from $Y$}.
}\end{Def}

Lemma~\ref{unbndlem} can be restated as saying that $\mathbf{H_1}$ is an unbounded distance from $\mathbf{H_2}$. 

\begin{Def}\textup{
Given $\delta \in \mathbf{R}^+$, a metric space $M$ is \textit{$\delta$-hyperbolic} if for any points $x_1, x_2, x_3 \in M$ and geodesics $x_1x_2, x_1x_3, x_2x_3$ between them, the geodesic $x_1x_2$ is contained in a $\delta$-neighborhood of 
$x_1x_3 \cup x_2x_3$.
}\end{Def}

\begin{Thm}[Masur and Minsky~\cite{mm:hyp}]
\label{mm1thm}
For each genus $g \geq 2$, there is a $\delta$ such that if $\Sigma$ is a closed genus-$g$ surface then $C(\Sigma)$ is $\delta$-hyperbolic.
\end{Thm}

\begin{Def}\textup{
A subset $A$ of $M$ is \textit{$T$-quasi-convex} if there is a $T \in \mathbf{R}^+$ such that for any points $x_1, x_2 \in A$, any geodesic $x_1x_2$ between the points is contained in a $T$-neighborhood of $A$.
}\end{Def}

\begin{Thm}[Masur and Minsky~\cite{mm:convex}]
\label{mm2thm}
For each $g \geq 2$, there is a $T$ such that for any $\phi : \Sigma \rightarrow H$, where $\Sigma$ is a closed genus-$g$ surface, the handlebody set corresponding to $H$ is $T$-quasi-convex.
\end{Thm} 

\begin{Def}\textup{
A subset $A$ of $M$ has \textit{infinite diameter} if for every $N \in \mathbf{N}$, there are points $x_1, x_2 \in A$ such that $d(x_1,x_2) > N$.
}\end{Def}

\begin{Thm}[Hempel~\cite{hemp:cc}]
\label{hemp1thm}
The handlebody set corresponding to a handlebody $H$ with genus $g \geq 2$ has infinite diameter.
\end{Thm}

\begin{Thm}[Hempel~\cite{hemp:cc}]
\label{hemp2thm}
For every $g \geq 2$ and for every $N$, there is a manifold with a genus-$g$ Heegaard splitting $(\Sigma', H_1', H_2')$ such that $d(\mathbf{H_1'}, \mathbf{H_2'}) > N$.
\end{Thm}

Consider the group $G$ of automorphisms of $\Sigma$.  An element $g \in G$ sends each isotopy class of curves in $\Sigma$ to a new isotopy class.  If two curves are disjoint then their images in $g$ will be disjoint, so $g$ defines an automorphism of $C(\Sigma)$.  

Given a handlebody set $\mathbf{H} \subset C(\Sigma)$, we will say that 
$g$ \textit{preserves} $\mathbf{H}$ if it sends every point of 
$\mathbf{H}$ to another point of $\mathbf{H}$.  Note that if $g$ is induced
by an automorphism of $H$ then $g$ will preserve $\mathbf{H}$.

Again, let $(\Sigma, H_1, H_2)$ be a Heegaard splitting of $S^3$ and let
$\mathbf{H_1}$, $\mathbf{H_2}$ be the handlebody sets corresponding to $H_1$,
$H_2$.  We will need the following Lemma:

\begin{Lem}
\label{genlem}
The group $G$ is generated by the automorphisms $\{g \in G : g(\mathbf{H_1}) = \mathbf{H_1} \ or\ g(\mathbf{H_2}) = \mathbf{H_2} \}$ (the automorphisms that perserve either $\mathbf{H_1}$ or $\mathbf{H_2}$).
\end{Lem}

\begin{figure}[htb]
  \begin{center}
  \includegraphics[width=3.5in]{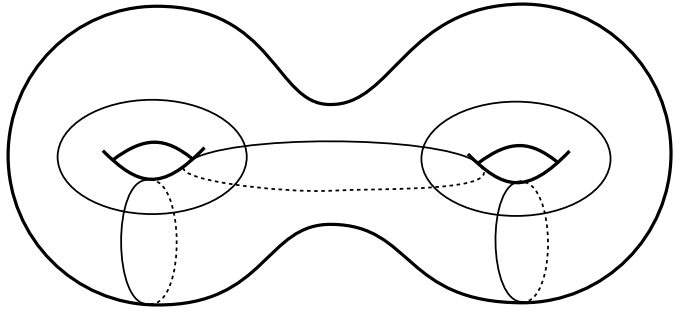}
  \caption{Dehn twists along the loops shown generate the mapping class group of a genus-two surface.  There is a similar collection for any higher genus surface as well.}
  \label{mcgfig}
  \end{center}
\end{figure}

\begin{proof}
The automorphism group of $\Sigma$ is generated by Dehn twists along a certain finite collection of simple closed curves curves. The loops for a genus-two surface are shown in Figure~\ref{mcgfig}.  Each loop bounds a disk in either $H_1$ or $H_2$ so the automorphism given by a Dehn twist along one of these loops extends to an automorphism of $H_1$ or $H_2$.  This automorphism preserves either $\mathbf{H_1}$ or $\mathbf{H_2}$.
\end{proof}

\begin{proof}[Proof of Lemma~\ref{unbndlem}]
Assume for contradiction there is an integer $N$ such that for every
vertex $v \in \mathbf{H_1}$, $d(v, \mathbf{H_2}) < N$.  There is an
automorphism of $S^3$ which swaps $H_1$ and $H_2$, so it must also be true
that for every $v \in \mathbf{H_2}$, $d(v, \mathbf{H_1}) < N$.

Let $g_1 \in G$ preserve $\mathbf{H_1}$.  The set $g_1(\mathbf{H_2})$ is a handlebody set of $C(\Sigma)$.  Since $g_1$ is an isometry of $C(\Sigma)$, we know that for any $v \in g_1(\mathbf{H_2})$, $d(v, \mathbf{H_1}) < N$. Let $v' \in \mathbf{H_1}$ such that $d(v, v') < N$.  Then $d(v', \mathbf{H_2}) < N$ so $d(v, \mathbf{H_2}) < 2N$.  Generalizing this argument, we see that if some subset $X$ is a bounded distance from $\mathbf{H_1}$ and from $\mathbf{H_2}$ and $g_1$ preserves $\mathbf{H_1}$ or $\mathbf{H_2}$ then $g_1(X)$ is a bounded distance from $\mathbf{H_1}$ and from $\mathbf{H_2}$ (although the bound is larger.) Let $g = g_1 \dots g_m \in G$ where each $g_i$ preserves either $\mathbf{H_1}$ or $\mathbf{H_2}$.  By induction, $g(\mathbf{H_2})$ is a bounded distance from $\mathbf{H_1}$ and from $\mathbf{H_2}$.

For any handlebody set $\mathbf{H}$ corresponding to a handlebody $H$ and a
homeomorphism $\phi : \Sigma \rightarrow \partial H$, there is a homeomorphism 
$\alpha : H \rightarrow H_2$ since $H$ and $H_2$ are handlebodies of the
same genus.  Let $i : \Sigma \rightarrow H_2$ be the inclusion map.  Then
$g = \phi^{-1} \circ \alpha \circ i$ is an automorphism of $\Sigma$
and $\mathbf{H} = g(\mathbf{H_2})$.

By Lemma~\ref{genlem}, $g$ is the composition of elements of $G$ which
preserve either $\mathbf{H_1}$ or $\mathbf{H_2}$ so the distance from $\mathbf{H}$ to $\mathbf{H_1}$ is bounded.  Since $\mathbf{H}$ was arbitrary, this implies that for any handlebody set in $C(\Sigma)$, the distance to $\mathbf{H_1}$ is bounded.  All that is needed to find a contradiction is to find some handlebody set $\mathbf{H}$ such that for every $N$ there is a $v \in \mathbf{H}$ with $d(v, H_1) > N$.

By Theorem~\ref{mm1thm} and Theorem~\ref{mm2thm}, let $\delta$ and $T$ be
integers such that $C(\Sigma)$ is $\delta$-hyperbolic and every handlebody
set is $T$-quasi-convex.  By Theorem~\ref{hemp2thm}, let  $M$ be a manifold with a Heegaard splitting $(\Sigma', H_1', H_2')$ such that $d(\mathbf{H_1'}, \mathbf{H_2'}) > 2 \delta + 2T$.

Identifying the handlebodies $H_1$ and $H_1'$ induces a homeomorphism of
from $\Sigma$ to $\Sigma'$ and therefore an isometry from $C(\Sigma')$ to
$C(\Sigma)$ which sends $\mathbf{H_1'}$ to $\mathbf{H_1}$.  The image
$\mathbf{H}$ of $\mathbf{H_2'}$ is such that 
$d(\mathbf{H}, \mathbf{H_1}) > 2 \delta + 2T$.

As we saw, the assumption that $\mathbf{H_1}$ is a bounded distance from
$\mathbf{H_2}$ implies that $\mathbf{H}$ is a bounded distance from
$\mathbf{H_1}$, i.e. there is an integer $R$ such that for
every $v \in \mathbf{H}$, $d(v, \mathbf{H_1}) < R$.  
By Theorem~\ref{hemp1thm}, let $v_1, v_2 \in \mathbf{H}$ such that
$d(v_1, v_2) > 2R + 4 \delta + 1$.  Let $v_3, v_4 \in \mathbf{H_1}$ such 
that $d(v_1, v_3) < R$ and $d(v_2, v_4) < R$.  

Let $v_1v_2$, $v_1v_3$, $v_1v_4$, $v_2v_4$ and $v_3v_4$ be geodesics, as in Figure~\ref{hyperectfig}. By Theorem~\ref{mm1thm}, $v_3v_4$ is contained in a $\delta$-neighborhood of $v_1v_3 \cup v_1v_4$ and $v_1v_4$ is contained in a $\delta$-neighborhood of $v_1v_2 \cup v_2v_4$.  Thus $v_3 v_4$ is contained in a $2 \delta$-neighborhood of $v_1v_3 \cup v_1v_2 \cup v_2v_4$.
\begin{figure}[htb]
  \begin{center}
  \includegraphics[width=2in]{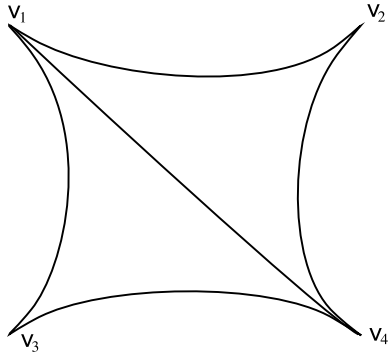}
  \caption{Geodesics connecting four vertices in the curve complex.}
  \label{hyperectfig}
  \end{center}
\end{figure}

If a vertex $x$ in $v_3v_4$ is in a $2 \delta$ neighborhood of $v_1v_2$ then there is a vertex $y$ in $v_1v_2$ such that $d(x, y) \leq 2 \delta$.   Because $\mathbf{H}$ and $\mathbf{H_1}$ are $T$-quasi-convex, this implies $d(\mathbf{H}, \mathbf{H_1}) \leq 2T + 2 \delta$.  By assumption, this inequality does not hold, so $v_3v_4$ is in a $2 \delta$-neighborhood of $v_1v_3 \cup v_2v_4$.

Because every vertex of the geodesic $v_3v_4$ is within $2 \delta$ of one of the geodesics, there must be adjacent vertices $x$ and $y$ such that $d(v_1v_3, x) \leq 2 \delta$ and $d(v_2v_4, y) \leq 2 \delta$.  Because each geodesic $v_1v_3$ and $v_2v_4$ has length at most $R$, there is a path from $v_1$ to $x$ with length at most $R + 2 \delta$ and a path from $y$ to $v_2$ of length at most $R + 2 \delta$.  Because $x$ and $y$ are adjacent, combining these paths produces a path from $v_1$ to $v_2$ of length at most $2R + 4 \delta + 1$.  This contradicts the assumption that $d(v_1,v_2) > 2R + 4 \delta + 1$ and the contradiction completes the proof. 
\end{proof}

Note that in the Heegaard splitting $(\Sigma', H_1', H_2')$ where $d(\mathbf{H_1'},\mathbf{H_2'}) > 2T + 2 \delta$, if we take a sequence of vertices $v_i \in \mathbf{H_2'}$ such that $d(v_1,v_i) \rightarrow \infty$ then we know that 
$d(v_i, \mathbf{H_1'}) \rightarrow \infty$.  However, for a Heegaard splitting $(\Sigma, H_1, H_2)$ of $S^3$, this in not necessarily true.  In fact the set $\mathbf{H_1} \cap \mathbf{H_2}$ appears to have infinite diameter.  Thus the proof does not suggest a way to actually contruct an example of a knot with $d(\tau) > N$.

\section{The Main Theorem}
\label{mainthmsect}

\begin{proof}[Proof of Lemma~\ref{highNlem}]
By Lemma~\ref{unbndlem}, there is a disk $D \subset H_1$ such that $d(\partial D, \mathbf{H_2}) > N + 1$.  If $D$ is non-separating, let $D_1 = D$.  Otherwise, let $D_1$ be a non-separating disk disjoint from $D$.  Such a disk exists and $d(\partial D_1, \mathbf{H_2}) > N$.  

Let $D_2, D_3$ be disjoint disks in $H_1$ which are non-separating and disjoint from $D_1$.  There is a spine of $H_1$ dual to the disks with edges $e_1, e_2, e_3$.  
Let $K = e_2 \cup e_3$ and $\tau = e_1$.  By definition $d(\tau) = d(\partial D_1, \mathbf{H_2}) > N$.
\end{proof}

For the proof of Theorem~\ref{mainthm}, we need one more Lemma:

\begin{Lem}
\label{torlem}
If the complement of $K$ is toroidal then $d(\tau) \leq 3$.
\end{Lem}

\begin{proof}
Let $D_1 \subset H_1$ be the disk dual to $\tau$.  Thompson~\cite{thm:dcp} showed that if $(\Sigma, H_1, H_2)$ is a Heegaard splitting of a toroidal manifold then $d(\Sigma) \leq 2$.  By Lemma~\ref{splitlem}, $d(\tau) \leq d(\Sigma) + 1 \leq 3$.
\end{proof}

In fact, toroidal, tunnel number one knots have been classified by Morimoto and Sakuma~\cite{ms:torknots}.  All such knots have $b_1(K) \leq 1$, implying $d(\tau) = 2$ for some tunnel.

\begin{proof}[Proof of Theorem~\ref{mainthm}]
Assume $N > 5$.  By Lemma~\ref{highNlem}, there is a knot $K$ with unknotting tunnel $\tau$ such that $d(\tau) > N$.  By Lemma~\ref{bridgelem}, $b_0(K) > N$ and by Lemma~\ref{torbridgelem}, $b_1(K) > 1$.

All that remains is to show that $K$ is hyperbolic.  If it's not, then either $K$ is a torus knot or its complement is toroidal.  Every torus knot has a one-bridge presentation with respect to an unknotted torus so $b_1(K) > 1$ implies $K$ is not a torus knot.  By Lemma~\ref{torlem}, the complement of $K$ will be atoroidal for $d(\tau) > 3$.  Since we assumed $N > 5$, we know that $K$ is hyperbolic.  This completes the proof.
\end{proof}

\bibliographystyle{abbrv}
\bibliography{bridges}

\end{document}